\documentclass{amsart}
\usepackage[latin1]{inputenc}
\usepackage{amsmath}
\usepackage{amsthm}
\usepackage{amsfonts}
\usepackage{amssymb}
\usepackage[T1]{fontenc}
\usepackage{float}
\usepackage[all]{xy}

\newtheorem*{theoremnn}{Theorem}
\newtheorem*{theoremmain}{Theorem \ref{theorem:main}}

\newtheorem{theorem}{Theorem}
\newtheorem{prop}[theorem]{Proposition}
\newtheorem{lem}[theorem]{Lemma}

\newtheorem{conj}[theorem]{Conjecture}

\theoremstyle{definition}
\newtheorem{defi}[theorem]{Definition}

\theoremstyle{definition}
\newtheorem{exmp}[theorem]{Example}

\def\Ext{{\rm{Ext}}}

\def\Gr{{\rm{Gr}}}

\def\End{{\rm{End}}}
\def\dim{{\rm{dim}\,}}
\def\ddim{{\textbf{dim}\,}}

\def\rep{{\rm{rep}}}
\def\ind{{\textrm{ind-}}}
\def\Ob{{\rm{Ob}}}

\def\N{{\mathbb{N}}}
\def\Z{{\mathbb{Z}}}
\def\C{{\mathbb{C}}}
\def\Q{{\mathbb{Q}}}

\def\CC{{\mathcal{C}}}

\def\<{\left<}
\def\>{\right>}

\def\ens#1{\left\{ #1 \right\}}
\def\fl{{\longrightarrow}\,}

\title{Positivity in coefficient-free rank two cluster algebras}
\author{\textsc{G. Dupont}}
 \address{Universit\'e de Lyon \\
Universit\'e Lyon 1 \\
Institut Camille Jordan CNRS UMR 5208 \\
43, boulevard du 11 novembre 1918\\
F-69622 Villeurbanne Cedex.}
 \email{dupont@math.univ-lyon1.fr}

\begin{document}
\maketitle

\begin{abstract}
	Let $b,c$ be positive integers, $x_1,x_2$ be indeterminates over $\Z$ and $x_m, m \in \Z$ be rational functions defined by $x_{m-1}x_{m+1}=x_m^b+1$ if $m$ is odd and $x_{m-1}x_{m+1}=x_m^c+1$ if $m$ is even. In this short note, we prove that for any $m,k \in \Z$, $x_k$ can be expressed as a substraction-free Laurent polynomial in $\Z[x_m^{\pm 1},x_{m+1}^{\pm 1}]$. This proves Fomin-Zelevinsky's positivity conjecture for coefficient-free rank two cluster algebras.
\end{abstract}

\setcounter{tocdepth}{1}
\tableofcontents

\section*{Introduction}
	\subsection*{A combinatorial result}
		Let $b,c$ be positive integers, the (coefficient-free) cluster algebra $\mathcal A(b,c)$ is a subring of the field $\Q(x_1,x_2)$ generated by the elements $x_m, m \in \Z$ satisfying the recurrence relations:
		$$x_{m+1}=\left\{\begin{array}{ll}
			\displaystyle \frac{x_m^b+1}{x_{m-1}} & \textrm{ if } m \in 2\Z+1~;\\
			\displaystyle \frac{x_m^c+1}{x_{m-1}} & \textrm{ if } m \in 2\Z.
		\end{array}\right.$$
		
		The elements $x_m, m \in \Z$ are called the \emph{cluster variables} of $\mathcal A(b,c)$ and the pairs $(x_m,x_{m+1}), m \in \Z$ are called the \emph{clusters} of $\mathcal A(b,c)$.
		
		The Laurent phenomenon \cite{cluster1} implies that for any $m \in \Z$ and any $k \in \Z$ the cluster variable $x_k$ belongs to the ring of Laurent polynomials $\Z[x_m^{\pm 1}, x_{m+1}^{\pm 1}]$. 
		
		When $bc \leq 4$, it was proved by Sherman-Zelevinsky \cite{shermanz} and independently by Musiker-Propp \cite{PM:rangdeux} that for any $m \in \Z$ and any $k \in \Z$ the cluster variable $x_k$ belongs to $\N[x_m^{\pm 1}, x_{m+1}^{\pm 1}]$. This was later proved by Caldero-Reineke for any $b=c$ \cite{CR}. In this paper, we prove this for arbitrary positive integers $b,c$. More precisely, the main result of the paper is :
		\begin{theoremmain}
			Let $b,c$ be positive integers. With the above notations, we have
			$$\ens{x_k \ : \ k \in \Z} \subset \N[x_m^{\pm 1}, x_{m+1}^{\pm 1}]$$
			for any $m \in \Z$.
		\end{theoremmain}
	
	\subsection*{The positivity conjecture for cluster algebras}
		In particular, this result is a particular case of a general conjecture formulated by Fomin and Zelevinsky for arbitrary cluster algebras. We recall that cluster algebras were introduced by Fomin and Zelevinsky in a series of papers \cite{cluster1,cluster2,cluster3,cluster4} in order to design a general framework for understanding total positivity in algebraic groups and canonical bases in quantum groups. They turned out to be related to various subjects in mathematics like combinatorics, Lie theory, representation theory, Teichm\"uller theory and many other topics. 
	
		In full generality, a (coefficient-free) cluster algebra is a commutative algebra generated by indeterminates over $\Z$ called \emph{cluster variables}. They are gathered into sets of fixed cardinality called \emph{clusters}. The initial data for constructing a (coefficient-free) cluster algebra is a \emph{seed}, that is, a pair $(B,\textbf u)$ where $B \in M_q(\Z)$ is a skew-symmetrizable matrix and $\textbf u=(u_1, \ldots, u_q)$ is a $q$-tuple of indeterminates over $\Z$. The cluster variables are defined inductively by a process called \emph{mutation}. The cluster algebra associated to a seed $(B,\textbf u)$ is denoted by $\mathcal A(B)$. For every cluster $\textbf c=\ens{c_1, \ldots, c_q}$ in $\mathcal A(B)$, the \emph{Laurent phenomenon} ensures that the cluster algebra $\mathcal A(B)$ is a $\Z$-subalgebra of the ring $\Z[c_1^{\pm 1}, \ldots, c_q^{\pm 1}]$ of Laurent polynomials in $\textbf c$ \cite{cluster1}. 
		
		A Laurent polynomial is called \emph{substraction-free} it it can be written as a $\N$-linear combination of unitary Laurent monomials. In \cite{cluster1}, Fomin and Zelevinsky gave the so-called \emph{positivity conjecture} for arbitrary cluster algebras. In the coefficient-free case, this conjecture can be stated as follows:
		\begin{conj}[\cite{cluster1}]\label{conj:positivity}
			Let $\mathcal A$ be a cluster algebra. Then any cluster variable $x$ in $\mathcal A$ can be written as a substraction-free Laurent polynomial in any cluster of $\mathcal A$.
		\end{conj}
		
		A simple but non-trivial class of cluster algebras is constituted by the so-called \emph{rank two cluster algebras}, that is, the cluster algebras of the form $\mathcal A(B_{b,c})$ with
		$$B_{b,c}=\left[\begin{array}{rr} 0 & b \\ -c & 0 \end{array}\right] \in M_2(\Z).$$
		Note that $\mathcal A(B_{b,c})$ is the cluster algebra $\mathcal A(b,c)$ introduced at the beginning.
		
		Thus, Theorem \ref{theorem:main} is equivalent to the following statement:
		\begin{theoremnn}
		 	Let $b,c$ be positive integers. Then Conjecture \ref{conj:positivity} holds for the cluster algebra $\mathcal A(b,c)$.
		\end{theoremnn}

	\subsection*{Organization of the paper}
		Despite the fact the main theorem can be expressed in purely combinatorial terms, the methods we use in this paper are based on representation theory and more precisely on categorifications of cluster algebras using cluster categories and cluster characters developed in \cite{BMRRT,CC,CK1,CK2}. Positivity will follow from results of Caldero-Reineke \cite{CR} and folding processes inspired by methods in \cite{Dupont:nonsimplylaced} (see also \cite{Demonet, Demonet:PhD}). In section \ref{section:clustercat}, we recall the necessary background on cluster categories and cluster characters. In section \ref{section:automorphisms}, we investigate briefly cluster categories and cluster algebras associated to quivers with automorphisms. In section \ref{section:unfolding}, we define a folding process in order to realize cluster variables in rank two cluster algebras and we prove the main result.

\begin{section}{Cluster categories and cluster characters}\label{section:clustercat}
	In this section, we recall necessary background on cluster categories and cluster characters. In the whole paper, $k$ will denote the field $\C$ of complex numbers. Let $Q=(Q_0,Q_1)$ be a finite acyclic quiver where $Q_0$ is the set of vertices and $Q_1$ the set of arrows. For an arrow $\alpha: i \fl j \in Q_1$, we denote by $s(\alpha)=i$ its \emph{source} and $t(\alpha)=j$ its \emph{target}. We assume that $Q$ has no oriented cycles, $Q$ is thus called \emph{acyclic}.
	
	To the quiver $Q$, we can associate a skew-symmetric matrix $B_Q=(b_{ij})_{i,j \in Q_0}$ as follows:
	$$b_{ij}=|\ens{\alpha \in Q_1 \ : \ s(\alpha)=i \textrm{ and } t(\alpha)=j}|-|\ens{\alpha \in Q_1 \ : \ t(\alpha)=i \textrm{ and } s(\alpha)=j}|.$$
	This induces a 1-1 correspondence from the set of quivers without loops and 2-cycles to the set of skew-symmetric matrices.
	
	\begin{subsection}{Cluster categories}
		A \emph{representation} $V$ of $Q$ is a pair $$V=((V(i))_{i \in Q_0},(V(\alpha))_{\alpha \in Q_1})$$ where $(V(i))_{i \in Q_0}$ is a family of finite dimensional $k$-vector spaces and $(V({\alpha}))_{\alpha \in Q_1}$ is a family of $k$-linear maps $V(\alpha):V(s(\alpha)) \fl V(t(\alpha))$. 
		A \emph{morphism of representations} $f:V \fl W$ is a family $(f_i)_{i \in Q_0}$ of $k$-linear maps such that the following diagram commutes
		$$\xymatrix{
			V(s(\alpha)) \ar[r]^{V(\alpha)} \ar[d]^{f_{s(\alpha)}}& V(t(\alpha)) \ar[d]^{f_{t(\alpha)}}\\
			W(s(\alpha)) \ar[r]^{W(\alpha)} & W(t(\alpha))
		}$$
		for any arrow $\alpha \in Q_1$. This defines a category $\rep(Q)$ which is equivalent to the category $kQ$-mod of finitely generated modules over the path algebra $kQ$ of $Q$. 
		
		For any vertex $i \in Q_0$, we denote by $P_i$ (resp. $I_i$, $S_i$) the indecomposable projective (resp. injective, simple) $kQ$ module associated to the vertex $i$. For any representation $M$ of $Q$, the \emph{dimension vector} of $M$ is the vector $\ddim M=(\dim M(i))_{i \in Q_0} \in \N^{Q_0}$.
		
		We denote by $D^b(kQ)$ the bounded derived category of $Q$. This is a triangulated category with shift functor $[1]$ and Auslander-Reiten translation $\tau$. The \emph{cluster category of $Q$}, introduced in \cite{BMRRT}, is the orbit category $\CC_Q=D^b(kQ)/\tau^{-1}[1]$ of the auto-functor $\tau^{-1}[1]$ in $D^b(kQ)$. This is a triangulated category \cite{K}. Moreover, it is proved in \cite{BMRRT} that indecomposable objects are given by
		$$\ind(\CC_Q)=\ind(kQ\textrm{-mod}) \sqcup \ens{P_i[1] \ : \ i \in Q_0}$$
		and that $\CC_Q$ is a 2-Calabi-Yau category, that is, that there is a functorial duality
		$$\Ext^1_{\CC_Q}(M,N) \simeq D \Ext^1_{\CC_Q}(N,M)$$
		for any two objects $M,N$ in the cluster category.
	\end{subsection}
	
	\begin{subsection}{The Caldero-Chapoton map}
		We denote by $\<-,-\>$ the Euler form on $kQ$-mod. For any representation $M \in \rep(Q)$ and any $\textbf e \in \N^{Q_0}$, the \emph{grassmannian of submodules} of $M$ is the projective variety
		$$\Gr_{\textbf e}(M)=\ens{N \textrm{ subrepresentation of M s.t. } \ddim N = \textbf e}.$$
		We can thus consider the Euler-Poincar\'e characteristic $\chi(\Gr_{\textbf e}(M))$ of this variety.
		
		We denote by $\mathcal A(Q)$ the coefficient-free cluster algebra with initial seed $(B_Q,\textbf u)$ where $\textbf u=\ens{u_i \ : \ i \in Q_0}$ is a set of indeterminates over $\Q$. In \cite{CC}, the authors considered a map $X_?:\Ob(\CC_Q) \fl \Z[\textbf u^{\pm 1}]$ which is now referred to as the \emph{Caldero-Chapoton map}.
		
		\begin{defi}
			The \emph{Caldero-Chapoton map}\index{Caldero-Chapoton map} is the map $X_?$ defined from the set of objects in $\mathcal C_Q$ to the ring of Laurent polynomials in the indeterminates $\ens{u_i, i \in Q_0}$ by:
			\begin{enumerate}
				\item[a.] If $M$ is an indecomposable $kQ$-module, then
					\begin{equation}\label{CCmap}
						X_M = \sum_{\textbf e} \chi(\Gr_{\textbf e}(M)) \prod_{i \in Q_0} u_i^{-<\textbf e, \ddim S_i>-\<\ddim S_i, \ddim M - \textbf e\>};
					\end{equation}
				\item[b.] If $M=P_i[1]$ is the shift of the projective module associated to $i \in Q_0$, then $$X_M=u_i;$$
				\item[c.] For any two objects $M,N$ in $\mathcal C_Q$, 
					$$X_{M \oplus N}=X_MX_N.$$
			\end{enumerate}
		\end{defi}
		Note that equality (\ref{CCmap}) holds also for decomposable modules.

		One of the main motivations for introducing the Caldero-Chapoton map was:
		\begin{theorem}[\cite{CK2}]\label{theorem:correspondenceCK2}
			Let $Q$ be an acyclic quiver. Then $X_?$ induces a 1-1 correspondence from the set of indecomposable objects without self-extensions in $\CC_Q$ and cluster variables in $\mathcal A(Q)$.
		\end{theorem}
		
		Moreover, it gives the cluster algebra $\mathcal A(Q)$ a structure of Hall algebra of the cluster category $\CC_Q$. More precisely, Caldero and Keller proved the following:
		\begin{theorem}[\cite{CK2}]\label{theorem:mult}
			Let $Q$ be an acyclic quiver, $M,N$ be two objects in $\CC_Q$ such that $\Ext^1_{\CC_Q}(M,N) \simeq k$, then 
			$$X_MX_N=X_B+X_{B'}$$
			where $B$ and $B'$ are the unique objects such that there exists non-split triangles
			$$M \fl B \fl N \fl M[1] \textrm{ and } N \fl B' \fl M \fl N[1].$$
		\end{theorem}
		
		Caldero and Reineke later proved an important result towards about positivity:
		\begin{theorem}[\cite{CR}]\label{theorem:CR}
			Let $Q$ be an acyclic quiver, $M$ be an indecomposable module, then $\chi(\Gr_{\textbf e}(M)) \geq 0$.
		\end{theorem}
		
		In particular, this result combined with Theorem \ref{theorem:correspondenceCK2} proves that if $Q$ is an acyclic quiver, cluster variables in $\mathcal A(Q)$ can be written as substraction-free expressions in the initial cluster.
	\end{subsection}
\end{section}

\begin{section}{Automorphisms of quivers}\label{section:automorphisms}
	Let $Q=(Q_0,Q_1)$ be an acyclic quiver and $B_Q =(b_{ij}) \in M_{Q_0}(\Z)$ be the associated matrix.
	
	A subgroup $G$ of the symmetric group $\mathfrak S_{Q_0}$ is called a \emph{group of automorphisms} of $Q$ if $b_{gi,gj}=b_{i,j}$ for any $i,j \in Q_0$ and any $g \in G$.
	
	\begin{subsection}{$G$-action on the cluster category}
		Fix an acyclic quiver $Q$ equipped with a group of automorphisms $G$. We define a group action of $G$ on $\rep(Q)$ as follows. For any $g \in G$ and any representation $V$, set $gV$ to be the representation $((V(g^{-1}i))_{i \in Q_0},(V(g^{-1}\alpha))_{\alpha \in Q_1})$. For any morphism of representation $f:V \fl W$, $gf$ is the morphism $gV \fl gW$ given by $gf=(f_{g^{-1}i})_{i \in Q_0}$. Each $g$ defines a $k$-linear auto-equivalence of $\rep(Q)$ with quasi-inverse $g^{-1}$. 
		
		Each $g$ induces an action by auto-equivalence on the bounded derived category $D^b(k Q)$ commuting with the shift functor $[1]$ and the Auslander-Reiten translation $\tau$. Thus, each $g \in G$ induces an action by auto-equivalence on the cluster category $\CC_Q = D^b(kQ)/\tau^{-1}[1]$. This action is additive and given on shift of projective objects by $gP_i[1] \simeq P_{gi}[1]$ for any $g \in G$ and $i \in Q_0$.
	\end{subsection}
	
	\begin{subsection}{$G$-action and variables}
		We still consider an acyclic quiver $Q$ equipped with a group of automorphisms $G$. Let $\mathcal A(Q)$ be the coefficient-free cluster algebra with initial seed $(B_Q,\textbf u)$ with $\textbf u=\ens{u_i \ : \ i \in Q_0}$. We define an action of $G$ by $\Z$-algebra homomorphisms on the ring of Laurent polynomials $\Z[\textbf u^{\pm 1}]$ by setting $gu_i=u_{gi}$. 
		
		The following lemma is straightforward (see \cite{Dupont:nonsimplylaced} for a proof).
		\begin{lem}\label{lem:xgm}
			Let $Q$ be an acyclic quiver equipped with a group $G$ of automorphisms. Then for any object $M$ in $\CC_Q$, we have
			$$gX_M=X_{gM}.$$
		\end{lem}
	\end{subsection}
\end{section}

\begin{section}{Unfolding rank two cluster algebras}\label{section:unfolding}
	The main idea of this paper comes from folding processes first developed by the author in \cite{Dupont:nonsimplylaced}. Fix $b,c$ two positive integers and $\textbf v=\ens{v_1, \ldots, v_b}$,  $\textbf w=\ens{w_1, \ldots, w_c}$ two finite sets. Let $K_{b,c}$ be the quiver given by
	$$(K_{b,c})_0=\textbf v \sqcup \textbf w$$
	and for every $i \in \ens{1, \ldots, b}$, $j \in \ens{1, \ldots, c}$, there is exactly one arrow $v_i \fl w_j$. The quiver $K_{b,c}$ can be represented as follows:
	$$
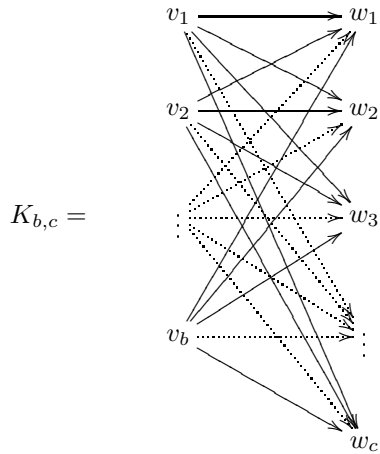
\xymatrix{
		&v_1 \ar[rr] \ar[rrd] \ar[rrdd] \ar@{.>}[rrddd] \ar[rrdddd] && w_1 \\
		&v_2 \ar[rru] \ar[rr] \ar[rrd] \ar@{.>}[rrdd] \ar[rrddd] && w_2 \\
		K_{b,c}=&\vdots \ar@{.>}[rruu] \ar@{.>}[rru] \ar@{.>}[rr] \ar@{.>}[rrd] \ar@{.>}[rrdd]&& w_3 \\
		&v_b \ar[rruuu] \ar[rruu] \ar[rru] \ar@{.>}[rr] \ar[rrd]&& \vdots \\
		& && w_c \\
	}$$
	
	Note that $K_{b,c}$ is acyclic and that the group
	$$G=\mathfrak S_{\textbf v} \times \mathfrak S_{\textbf w}$$
	is a group of automorphisms for $K_{b,c}$.
	
	In order to simplify notations, we write $Q=K_{b,c}$.  We denote by $\mathcal A(Q)$ the coefficient-free cluster algebra with initial seed $\left(\textbf u, B_Q\right)$ where $$\textbf u=\ens{u_{v_1}, \ldots, u_{v_b},u_{w_1}, \ldots, u_{w_c}}.$$
	
	We define the following $\Z$-algebra homomorphism $\pi$ called \emph{folding}:
	$$\pi:\left\{\begin{array}{rcll}
		\Z[u_i^{\pm 1} \ : \ i \in Q_0] & \fl & \Z[x_1^{\pm 1}, x_2^{\pm 1}]\\
		u_{v_i} & \mapsto & x_1 & \textrm{ for all }i=1, \ldots, b;\\
		u_{w_i} & \mapsto & x_2 & \textrm{ for all }i=1, \ldots, c.\\
	\end{array}\right.$$

	We will denote by $P_{\textbf v}$ (resp. $P_{\textbf w}$) a representative of the set $\ens{P_{v_1}, \ldots, P_{v_b}}$ (resp. $\ens{P_{w_1}, \ldots, P_{w_c}}$).
		
	We denote by $X_?$ the Caldero-Chapoton map on $\CC_Q$. For every $g \in G$ and $i \in \Z$, we have $gP_{\textbf v}[i] \in \ens{P_{v_1}[i], \ldots, P_{v_b}[i]}$ and $gP_{\textbf w}[i] \in \ens{P_{w_1}[i], \ldots, P_{w_c}[i]}$. It follows from Lemma \ref{lem:xgm} that $\pi(P_{\textbf v}[i])$ and $\pi(P_{\textbf w}[i])$ are well-defined elements in $\Z[x_1^{\pm 1}, x_2^{\pm 1}]$.
	
	We can thus give the following description of cluster variables in $\mathcal A(b,c)$:
	\begin{prop}\label{prop:rang2}
		Let $b,c$ be positive integers. Then, for any $m \in \Z$, we have
		$$\begin{array}{rl}
			x_{2m+1} &=\pi(X_{P_{\textbf v}[m+1]})\\
			x_{2m+2} &=\pi(X_{P_{\textbf w}[m+1]})\\
		\end{array}$$
	\end{prop}
	\begin{proof}
		We prove it by induction on $m$. We have $x_1=\pi(u_{v_i})=\pi(X_{P_{v_i}[1]})$ for every $i=1, \ldots, b$ and $x_2=\pi(u_{w_i})=\pi(X_{P_{w_i}[1]})$ for every $i=1, \ldots, c$. 
		
		Fix $i \in\ens{1, \ldots, c}$ and $m \in \Z$, we have isomorphisms of $k$-vector spaces
		\begin{align*}
			k 
				& \simeq \End_{\CC_Q}(P_{w_i}[m+1])\\
				& \simeq \Ext^1_{\CC_Q}(P_{w_i}[m+1], P_{w_i}[m])\\
				& \simeq \Ext^1_{\CC_Q}(P_{w_i}[m], P_{w_i}[m+1])\\
		\end{align*}
		and the corresponding triangles are
		$$P_{w_i}[m] \fl 0 \fl P_{w_i}[m+1] \fl P_{w_i}[m+1],$$
		$$P_{w_i}[m+1] \fl \bigoplus_{j=1}^c P_{v_j}[m] \fl P_{w_i}[m] \fl P_{w_i}[m+2].$$
		
		It thus follows from Theorem \ref{theorem:mult} that
		$$X_{P_{w_i}[m]}X_{P_{w_i}[m+1]}=\prod_{i=1}^c X_{P_{v_j}[m]}+1$$
		but $\pi(X_{P_{w_i}[m]})=x_{2m}$ and $\pi(X_{P_{v_j}[m]})=x_{2m-1}$ for every $j=1,\ldots,b$. We thus get
		$$\pi(X_{P_{w_i}[m+1]})=\frac{x_{2m-1}^b+1}{x_{2m}}=x_{2m+2}.$$
		The other cases are proved similarly.
	\end{proof}
	
	\begin{figure}[htb]
		\setlength{\unitlength}{.3mm}
		\begin{picture}(300,200)(-80,-100)
			\multiput(-120,0)(120,0){3}{
				\put(0,60){\line(1,-1){50}}
				\put(0,60){\line(5,-6){50}}
				\put(50,0){\line(0,1){10}}
				
				\put(60,10){\line(1,1){50}}
				\put(60,0){\line(5,6){50}}
				\put(60,0){\line(0,1){10}}
			}
			\put(-130,65){$P_{\textbf v}[2]$}
			\put(-70,-7){$P_{\textbf w}[1]$}	
			\put(-10,65){$P_{\textbf v}[1]$}
			\put(50,-7){$P_{\textbf w}$}	
			\put(110,65){$P_{\textbf v}$}
			\put(170,-7){$P_{\textbf w}[-1]$}	
			\put(230,65){$P_{\textbf v}[-1]$}
			
			\put(270,50){$\Gamma(kQ)$}
			
			\put(-135,0){$\pi \circ X_?$}
			\put(-125,-40){\line(0,1){85}}
			\put(-65,-40){\line(0,1){27}}	
			\put(-5,-40){\line(0,1){85}}
			\put(45,-40){\line(0,1){27}}	
			\put(115,-40){\line(0,1){85}}
			\put(175,-40){\line(0,1){27}}
			\put(235,-40){\line(0,1){85}}

			\put(-130,-50){$x_3$}
			\put(-70,-50){$x_2$}	
			\put(-10,-50){$x_1$}
			\put(50,-50){$x_0$}	
			\put(110,-50){$x_{-1}$}
			\put(170,-50){$x_{-2}$}	
			\put(230,-50){$x_{-3}$}
			\put(270,-50){$\mathcal A(b,c)$}
		
		\end{picture}
		\caption{Realizing cluster variables}\label{figure:quotientrangdeux}
	\end{figure}
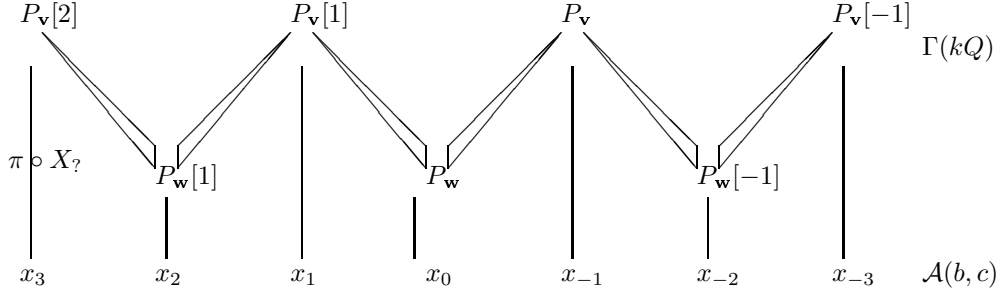
	
	Figure \ref{figure:quotientrangdeux} sums up the situation where in the AR-quiver of $\CC_Q$, we grouped together the objects in a same $G$-orbit and the arrows between these orbits. 
	
	Locally, the quiver has actually the following shape:
	$$
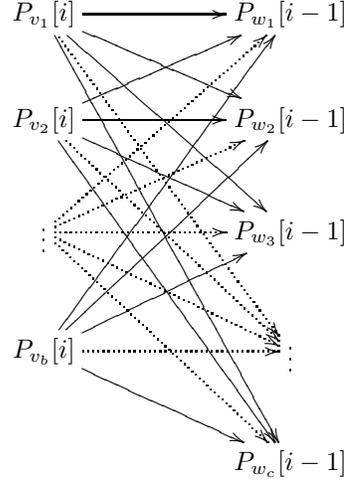
\xymatrix{
		&P_{v_1}[i] \ar[rr] \ar[rrd] \ar[rrdd] \ar@{.>}[rrddd] \ar[rrdddd] && P_{w_1}[i-1] \\
		&P_{v_2}[i] \ar[rru] \ar[rr] \ar[rrd] \ar@{.>}[rrdd] \ar[rrddd] && P_{w_2}[i-1] \\
		&\vdots \ar@{.>}[rruu] \ar@{.>}[rru] \ar@{.>}[rr] \ar@{.>}[rrd] \ar@{.>}[rrdd]&& P_{w_3}[i-1] \\
		&P_{v_b}[i] \ar[rruuu] \ar[rruu] \ar[rru] \ar@{.>}[rr] \ar[rrd]&& \vdots \\
		& && P_{w_c}[i-1] \\
	}$$
	
	We now prove the main theorem:
	\begin{theorem}\label{theorem:main}
		Let $b,c$ be positive integers. With the above notations, we have
		$$\ens{x_k \ : \ k \in \Z} \subset \N[x_m^{\pm 1}, x_{m+1}^{\pm 1}]$$
		for any $m \in \Z$.
	\end{theorem}
	\begin{proof}
		 Let $b,c$ be positive integers. All the exchange matrices in $\mathcal A(b,c)$ are either $\left[\begin{array}{rr} 0 & b \\ -c & 0 \end{array}\right]$ or $\left[\begin{array}{rr} 0 & -b \\ c & 0 \end{array}\right]$. It is thus enough to prove that cluster variables in $\mathcal A(b,c)$ can be expressed as substraction-free expressions in the initial cluster $\textbf c=(x_1,x_2)$.
		 
		Fix a cluster variable $x$ in $\mathcal A(b,c)$. Then $x=x_k$ for some $k \in \Z$. Assume for example that $k=2m+2$ for some $m \in \Z$. Then it follows from Proposition \ref{prop:rang2} that $x=\pi(X_{P_{\textbf w}[m+1]})$. According to Theorem \ref{theorem:CR}, as $K_{b,c}$ is acyclic, $X_{P_{\textbf w}[m+1]}$ is a linear combination of Laurent monomials in $\textbf u$ with positive coefficients. Thus, applying $\pi$ to this expansion, $x=\pi(X_{P_{\textbf w}[m+1]})$ is a linear combination of Laurent monomials in $(x_1,x_2)$ with positive coefficients. Similarly, if $k=2m+1$, $x_k$ can be written as a linear combination of Laurent monomials in $(x_1,x_2)$ with positive coefficients. This proves the theorem.
	\end{proof}
	
	\begin{exmp}
		We consider the cluster algebra $\mathcal A(2,3)$ with initial cluster $(x_1,x_2)$. We consider the following quiver:
		$$
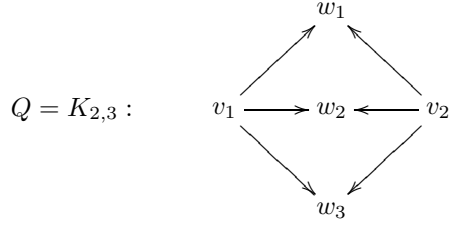
\xymatrix{
				&& w_1 \\
			Q=K_{2,3}: & v_1 \ar[ru]\ar[r]\ar[rd] & w_2 & \ar[l]\ar[ld]\ar[lu] v_2\\
				&& w_3
		}$$
		
		Let $\mathcal A(Q)$ be the cluster algebra with initial seed $(B_Q,\textbf u)$ where
		$$\textbf u=\ens{u_{v_1}, u_{v_2}, u_{w_1}, u_{w_2}, u_{w_3}}.$$ 
		Let $X_?: \Ob(\CC_Q) \fl \Z[\textbf u^{\pm 1}]$ be the corresponding Caldero-Chapoton map.
		
		A direct computation gives
		$$X_{P_{w_j}}=\frac{1+u_{v_1}u_{v_2}}{u_{w_j}},$$
		$$X_{I_{w_j}}=X_{P_{w_j}[2]}=\frac{1+u_{v_1}u_{v_2}+2u_{w_1}u_{w_2}u_{w_3}+u_{w_1}^2u_{w_2}^2u_{w_3}^2}{u_{w_j}u_{v_1}u_{v_2}}$$
		for every $j=1,2,3$ and
		$$X_{P_{v_i}}=\frac{1+u_{v_1}^3u_{v_2}^3+3u_{v_1}^2u_{v_2}^2+3u_{v_1}u_{v_2}+u_{w_1}u_{w_2}u_{w_3}}{u_{v_i}u_{w_1}u_{w_2}u_{w_3}},$$
		$$X_{I_{v_i}}=X_{P_{v_i}[2]}=\frac{1+u_{w_1}u_{w_2}u_{w_3}}{u_{v_i}}$$
		for every $i=1,2$.
		
		Consider the folding morphism 
		$$\pi:\left\{\begin{array}{rcll}
			\Z[u_i^{\pm 1} \ : \ i \in Q_0] & \fl & \Z[x_1^{\pm 1}, x_2^{\pm 1}]\\
			u_{v_i} & \mapsto & x_1 & \textrm{ for all }i=1,2;\\
			u_{w_j} & \mapsto & x_2 & \textrm{ for all }j=1, 2, 3.\\
		\end{array}\right.$$
		
		Thus we get
		$$\pi(X_{P_{v_i}})=\frac{1+x_1^6+3x_1^4+3x_1^2+x_2^3}{x_1x_2^3}=x_{-1}, \quad \pi(X_{P_{w_j}})=\frac{1+x_1^2}{x_2}=x_0,$$
		$$\pi(X_{P_{v_i}[1]})=x_1, \quad \pi(X_{P_{w_j}[1]})=x_2,$$
		$$\pi(X_{I_{v_i}})=\frac{1+x_2^3}{x_1}=x_{3}, \quad \pi(X_{I_{w_j}})=\frac{1+x_1^2+2x_2^3+x_2^6}{x_1^2x_2}=x_4.$$
		This illustrates Theorem \ref{theorem:main}.
	\end{exmp}
\end{section}


\newcommand{\etalchar}[1]{$^{#1}$}

\end{document}